\documentclass[
final
]{dmtcs-episciences}

\usepackage[utf8]{inputenc}
\usepackage{subfigure}
\usepackage{amsmath}
\usepackage{bm}
\usepackage{float}
\usepackage[linguistics]{forest}
\def\wasc{{\textsf{wasc}}}
\newcommand{\A}{{\mathcal{W}}}
\newcommand{\wa}{{\bm w}}
\newcommand{\T}{{\mathcal{T}}}
\newcommand{\I}{{\mathcal I}}
\newcommand{\newmatrices}{\mathrm{WMat}}

\newcommand{\seqnum}[1]{\href{http://oeis.org/#1}{\underline{#1}}}

\usepackage{amsthm}
\theoremstyle{plain}
\newtheorem{theorem}{Theorem}[section]
\newtheorem{proposition}[theorem]{Proposition}

\newtheorem{observation}[theorem]{Observation}
\newtheorem{conjecture}[theorem]{Conjecture}
\theoremstyle{definition}
\newtheorem{definition}[theorem]{Definition}
\newtheorem{example}[theorem]{Example}

\usepackage[round]{natbib}

\author{Be\'ata B\'enyi}
\address{\noindent Department of Hydraulic Engineering, University of Public Service, Baja, HUNGARY}
\email{benyi.beata@uni-nke.hu}

\author[Be\'ata B\'enyi et. al]{Be\'ata B\'enyi\affiliationmark{1}
  \and Toufik  Mansour\affiliationmark{2}
  \and Jos\'e L. Ram\'{\i}rez\affiliationmark{3}}
\title[Pattern Avoidance in Weak Ascent Sequences]{Pattern Avoidance in Weak Ascent Sequences}
\affiliation{
  Department of Hydraulic Engineering, University of Public Service, Baja, Hungary \\
Department of Mathematics, University of Haifa, Haifa, Israel\\
Departamento de Matem\'aticas,  Universidad Nacional de Colombia,  Bogot\'a, Colombia}
\begin{document}

\publicationdata{vol. 26:1, Permutation Patterns 2023}{2024}{2}{10.46298/dmtcs.12273}{2023-09-14; 2023-09-14; 2024-02-02; 2024-03-30}
{2024-04-02}
\maketitle
\begin{abstract}~ 

In this paper,  we study pattern avoidance in weak ascent sequences, giving some results for patterns of length 3.
This is an analogous study to one given by Duncan and Steingr\'imsson (2011) for ascent sequences. More precisely, we provide systematically the generating functions for the number of weak ascent sequences avoiding the patterns $001, 011, 012, 021$, and $102$. Additionally, we establish bijective connections between pattern-avoiding weak ascent sequences and other combinatorial objects, such as compositions, upper triangular 01-matrices, and plane trees.
\end{abstract}

\section{Introduction}

Weak ascent sequences were introduced recently by  \cite{BCD} as a family of integer sequences that contains the well-known class of ascent sequences introduced by \cite{Bous}. The weak ascent sequences encode  a collection of permutations avoiding a length-4 bivincular
pattern, upper-triangular binary matrices that satisfy a column-adjacency rule, and factorial posets that are weakly (3+1)-free.

A \emph{weak ascent sequence} is a word $x=x_1x_2\cdots x_n$ over the set of non-negative integers satisfying $x_1=0$ and $x_i\leq \wasc(x_1x_2\cdots x_{i-1})+1$ for $i=2, \dots, n$, where $\wasc(x_1x_2\cdots x_\ell)$ is the number of \emph{weak ascents} in the word  $x_1x_2\cdots x_\ell$, that is, the number of positions $j$ such that $x_j\leq x_{j+1}$.  We denote the set of all weak ascent sequences of length $n$ by $\A_n$. By definition we set $\A_0=\{\epsilon\}$, where $\epsilon$ is the empty word,  the unique word of length zero. For example, the weak ascent sequences of length 4 are
\begin{multline*}
\A_4=\{0000, 0001, 0002, 0003, 0010, 0011, 0012, 0013, 0020, 0021,   0022, 0023,  \\ 0100, 0101, 0102, 0110, 0111, 0112, 0113, 0120, 0121, 0122, 0123\}.
\end{multline*}
A \emph{pattern} $p$ is a word containing each of the letters in $\{0,1, \dots , i-1\}$ at least once for some integer $i\geq 1$.  This is analogous to the classical definition of patterns for permutations, but here patterns may contain repeated letters.  A weak ascent sequence $x=x_1x_2\cdots x_n\in \A_n$ \emph{contains} a pattern $p=p_1\cdots p_k $  if there is a subsequence $x_{i_1}x_{i_2}\cdots x_{i_k}$ of $x$, $i_1<i_2<\cdots<i_k$, which is order-isomorphic to $p$ (appear in the same order of size). We say that $x$ \emph{avoids} the pattern $p$ if $x$ does not contain any occurrence of $p$. For example, the weak ascent sequence $0120244\in \A_{7}$   contains five occurrences of the pattern 001 and avoids  the pattern 312.  Let  $\A_n(p)$ denote the set of weak ascent sequences of length $n$ that avoid the pattern $p$ and  let $\wa_{p}(n)$ be the cardinality of  $\A_n(p)$. Further, we set $\A(p)=\bigcup_{n\geq 0}\A_n(p)$.

The goal of this work is to study weak ascent sequences avoiding certain patterns of length three. This study is analogous to the one conducted by \cite{DS} on ascent sequences. Our goal in this paper is to derive the generating function of some cases and show some connections with other combinatorial objects. Our study is far from comprehensive on the topic.

The remainder of this paper is structured as follows. In Section 2,  we provide enumerative results for the number of weak ascent sequences  avoiding the patterns $001, 011, 012, 021$, and $102$.  These enumerations correspond to known sequences in The On-Line Encyclopedia of Integer Sequences  \cite{OEIS}. Additionally,  we establish  bijections  with other combinatorial structures. In Section 3, we present our computer calculations concerning other patterns  (see Table \ref{table2})  and propose a conjecture regarding the case $210$, which proves to be more challenging than it appears at first glance.

\section{Enumerations}

The goal of the current section is to enumerate weak ascent sequences avoiding each given pattern in Table \ref{table1}.

\begin{table}[ht]
\centering
\resizebox{0.98\textwidth}{!}{
\begin{tabular}{|c|c|c|c|} \hline
Pattern $p$ & Number of non-empty sequences avoiding $p$ & OEIS & Formula/Generating function\\ \hline \hline
001 & 1, 2, 4, 8, 16, 32, 64, 128, 256, 512, $\dots$ & \seqnum{A000079} & $2^{n-1}$  \\ \hline

011 &  1, 2, 5, 14, 43, 143, 510, 1936, 7775, 32869& \seqnum{A098569} &$\sum_{m=1}^{n}\binom{n-1+\binom{m}{2}}{n-m}$ \\ \hline

012 & 1, 2, 5, 13, 34, 89, 233, 610, 1597, 4181,$\dots$ & \seqnum{A001519} & $F_{2n-1}$ (Fibonacci) \\ \hline

021 & 1, 2, 6, 21, 80, 322, 1347, 5798, 25512, 114236, $\dots$& \seqnum{A106228} & $ F_{021}(x) = 1 + \frac{xF_{021}(x)}{1 - xF_{021}(x)^2}$   \\ \hline

102 & 1, 2, 6, 22, 89, 381, 1694, 7744, 36168, 171831& \seqnum{A200753} & $F_{102}(x) = 1+ x(1-x)F_{102}(x)^3$.\\ \hline
\end{tabular}}
\caption{Weak ascent sequences avoiding a pattern of length three.}\label{table1}
\end{table}

\subsection{Pattern 001}

\begin{figure}[H]
\centering
\begin{forest}
for tree={fit=band,}
[0,[00,[000,[0000]]] [01,[010,[0100]] [011,[0110] [0111]] [012,[0120] [0121] [0122] [0123]]]]
\end{forest}
\caption{Generating tree $\T(\{001\})$.} \label{fig1}
\end{figure}
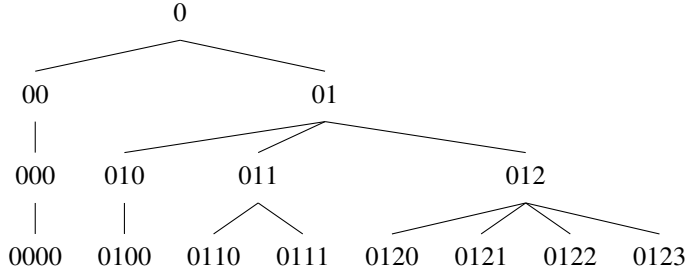

The simplest case is probably the pattern $001$.
We apply here the generating tree method introduced by  \cite{W}. We recall the definition for the general case. Given a fixed set of patterns, $B$, let $\A(B)$ denote the set of all non-empty weak ascent sequences that avoid every pattern contained in $B$. We define the generating tree (see \cite{W}) $\T(B)$ to be the plane tree as follows. We say $\T(B)$ is empty if there is no weak ascent sequence of arbitrary length avoiding the set $B$, that is, $0\in B$. Otherwise, the root can always be taken as $0$. Starting with this root which stays at level $1$, we construct the remainder of the nodes of the tree $\T(B)$ recursively as follows. The children of $x_1\cdots x_n\in\A_n(B)$ are those sequences from the set $$\{x_1\cdots x_{n+1}\mid x_{n+1}=0,1,\ldots,\wasc(x_1\cdots x_n)+1\}$$ that avoid the patterns in $B$.  For example, Figure \ref{fig1} shows the first levels of the tree $\T(\{001\})$.

We define an equivalence relation on the nodes of $\T(B)$. Let $\T(B;x)$ be the subtree consisting of the weak ascent sequence $x$ as the root and its descendants in $\T(B)$.  We say that $x$ is {\em equivalent} to $x'$ if and only if $\T(B;x)\cong\T(B;x')$ (in the sense of plane trees). Let $\T'(B)$ be the same tree as $\T(B)$, where we replace each node $x$ by the first node $x'\in\T(B)$ from top to bottom and from left to right in $\T(B)$ such that $\T(B;x)\cong\T(B;x')$. From now, we consider $\T'(B)$ and $\T(B)$ identical.

\begin{theorem}\label{th3}
The generating function for the number of non-empty weak ascent sequences of length $n$ that avoid the pattern $001$ is given by $x/(1-2x)$.
\end{theorem}
\begin{proof}
We consider the generating tree for the case $B=\{001\}$. The root of $\T(B)$ is $a_0=0$. The children of the root are $b_0=00$ and $a_1=01$, and their children are $000$ and $010, 011, 012$, respectively (see Figure \ref{fig1}). Since $\T(B;000)\cong\T(B;00)$ and $\T(B;010)\cong\T(B;00)$, we see that the children of the root $a_0$ in $\T(B)$ are $b_0,a_1$ and their children are $a_0$ and $b_0, b_1, a_2$, respectively. More generally, let $a_m=012\cdots m$  and $b_m=012\cdots mm$, then we see that $\T(B)$ satisfies the succession rules
$$a_m\rightsquigarrow b_0,b_1,\ldots,b_m,a_{m+1} \quad \mbox{and}
\quad b_m\rightsquigarrow b_0,b_1,\ldots,b_m.$$
Define $A_m(x)$ (respectively, $B_{m}(x)$) to be the generating function for the number of nodes at level $n\geq1$ for the subtree of $\T(B;a_m)$ (respectively, $\T(B;b_m)$), where its root stays at level $1$. Thus, the succession rules of $\T(B)$ give
$$A_m(x)=x+x\sum_{j=0}^mB_j(x)+xA_{m+1}(x)\quad \mbox{and}\quad
\quad B_m(x)=x+x\sum_{j=0}^mB_j(x).$$
Define $A(v)=\sum_{m\geq0}A_m(x)v^m$ and $B(v)=\sum_{m\geq0}B_m(x)v^m$. Then
$$A(v)=\frac{x}{1-v}+\frac{x}{1-v}B(v)+\frac{x}{v}(A(v)-A(0))
\quad \mbox{and} \quad B(v)=\frac{x}{1-v}+\frac{x}{1-v}B(v).$$
Therefore, $B(v)=x/(1-x-v)$, and the first equation with $v=x$ gives
\[A(0)=\frac{x}{1-x}(1+B(x))=\frac{x}{1-2x}.\qedhere\]
\end{proof}

From the generating function it is clear that $\wa_{001}(n)=2^{n-1}$ for all $n\geq 1$. For example, the weak ascent sequences corresponding to $n=4$ are the $8$ nodes at level $4$ of the generating tree $\T(\{001\})$ in Figure \ref{fig1}.  

\subsubsection{A map}
After obtaining this simple sequence for the numbers, we also describe a map from the set of $01$ bitsequences of length $n-1$ to the set of weak ascent sequences of length $n$ avoiding the pattern $001$.

Given a $01$ bitsequence of length $n-1$, associate first to it a path: for a $1$ an up step, for a $0$ a flat step. Then, label the nodes of the path with their height. Record the labels (the heights) now in the following order, first each label in the order as it occurs the first time, and then record the labels (that already occurred once) in non-increasing order.

A sequence obtained in this way is a weak ascent sequence,  as we can see that as follows. Let $P$ be the path, $j$ the height of the path, and $x$ the image of the map. First we show that $x\in \A_n$. The first $j+1$ elements of $x$ are $01\cdots j$, i.e., integers up to $j$ listed in increasing order. Clearly, \linebreak $x_{i}\leq \wasc(x_1\cdots x_{i-1})+1$, for $i\leq j+1$. Note, that $\wasc(x_1\cdots x_{j+1})=j$. Further, from the definition of the map it follows that for $i\geq j+2$ the values are necessarily $x_i\leq j$. Hence, for $i\geq j+2$ it is true that $x_{i}\leq \wasc(x_1\cdots x_{i-1})+1$, and we conclude $x\in \A_n$.

The pattern $001$ can not occur in $x_1\cdots x_{j+1}$, since it is a strictly increasing sequence of integers. Hence, if a pattern $001$ occurs in $x$, the second element (corresponding to the second zero in the pattern) has an index greater than $j+1$. However, the values after the $(j+1)$th element are ordered non-increasingly, so we can not find an element that would correspond to $1$ in an occurrence of the pattern $001$. We see that $x\in \A_n(001)$.

The steps of the map can be inverted. Let $x$ be a weak ascent sequence avoiding $001$ with a maximal strictly increasing prefix of length $j+1$. Draw a path of height $j$ with so many flat steps at each level as many times the value of the level occurs in the sequence $x$ after the $(j+1)$th element. The path containing flat steps and up steps can easily be decoded back to a $01$ bitsequence.

\begin{figure}[H]
\setlength{\unitlength}{0.5cm}
\centering
\begin{picture}(4,3)
	\put(0,0){\line(1,0){1}}\put(1,0){\line(1,1){1}}\put(2,1){\line(1,0){1}}\put(3,1){\line(1,1){1}}
	\put(0,0){\circle*{0.2}}\put(1,0){\circle*{0.2}}\put(2,1){\circle*{0.2}}\put(3,1){\circle*{0.2}}\put(4,2){\circle*{0.2}}
	\put(0,0.1){$0$}\put(1,0.1){$0$}\put(2,1.1){$1$}\put(3,1.1){$1$}\put(4,2.1){$2$}
\end{picture}\textcolor{red}{$01210$}\quad
\begin{picture}(4,3)
	\put(0,0){\line(1,0){1}}\put(1,0){\line(1,1){2}}\put(3,2){\line(1,0){1}}
	\put(0,0){\circle*{0.2}}\put(1,0){\circle*{0.2}}\put(2,1){\circle*{0.2}}\put(3,2){\circle*{0.2}}\put(4,2){\circle*{0.2}}
	\put(0,0.1){$0$}\put(1,0.1){$0$}\put(2,1.1){$1$}\put(3,2.1){$2$}\put(4,2.1){$2$}
\end{picture}\textcolor{red}{$01220$}\quad
\begin{picture}(4,3)
	\put(0,0){\line(1,1){1}}\put(1,1){\line(1,0){2}}\put(3,1){\line(1,1){1}}
	\put(0,0){\circle*{0.2}}\put(1,1){\circle*{0.2}}\put(2,1){\circle*{0.2}}\put(3,1){\circle*{0.2}}\put(4,2){\circle*{0.2}}
	\put(0,0.1){$0$}\put(1,1.1){$1$}\put(2,1.1){$1$}\put(3,1.1){$1$}\put(4,2.1){$2$}
\end{picture}\textcolor{red}{$01211$}
\begin{picture}(4,3)
	\put(0,0){\line(1,1){1}}\put(1,1){\line(1,0){1}}\put(2,1){\line(1,1){1}}\put(3,2){\line(1,0){1}}
	\put(0,0){\circle*{0.2}}\put(1,1){\circle*{0.2}}\put(2,1){\circle*{0.2}}\put(3,2){\circle*{0.2}}\put(4,2){\circle*{0.2}}
	\put(0,0.1){$0$}\put(1,1.1){$1$}\put(2,1.1){$1$}\put(3,2.1){$2$}\put(4,2.1){$2$}
\end{picture}\textcolor{red}{$01221$}
\caption{Examples for the bijection.} \label{fig-001pattern}
\end{figure}
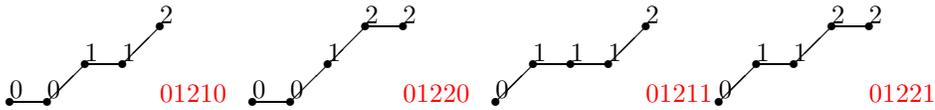

Figure \ref{fig-001pattern} shows some examples of the map for better understanding. Writing $F$ for a flat step and $U$ for an up step, in Figure \ref{fig-001pattern}
 the paths \textcolor{blue}{FUFU}, \textcolor{blue}{FUUF}, \textcolor{blue}{UFFU}, \textcolor{blue}{UFUF} and the corresponding weak ascent sequences are given.

After understanding the structure of the sequences in this class we provide also a direct counting. (This was pointed out by one of the referees.) Namely, let $j+1$ be the length of the maximal strictly increasing prefix of a weak ascent sequence $x$ of length $n$ avoiding $001$. Then the remainder of $x$ can consist precisely of any weakly decreasing sequence of the letters $\{0,1,\dots,j\}$, of which there are $\binom{n-1}{j}$. After summing up on $0\leq j\leq n-1$, we obtain $2^{n-1}$.

\subsection{Pattern 102}

Let $\I_n$ denote the set of \emph{inversion sequences} of length $n$, that is, integer sequences $e=e_1e_2\cdots e_n$ such that $0\leq e_i<i$ for each $i$.   \cite{CMSW} and  \cite{MS} began the study of classical patterns in inversion sequences. Denote by $\I_n(p)$ the set of inversion sequences of length $n$ that avoid $p$.  By the definitions we obtain the equality $\A_n(102)=\I_n(102)$ for all positive integer $n$. Indeed,  let $x$ be a sequence  in $ \I_n(102)$. Assume that $x$ has a position $j$ such that $\wasc(x_1\cdots x_{j-1})+1<x_j<j$. Then, the word $x_1\cdots x_{j-1}$ contains a descent, so  $x_1\cdots x_{j-1}x_{j}$ contains the pattern 102. This gives a contradiction, therefore $x\in \A_n(102)$. The other containment is clear and so the two sets are in fact equal.

From Theorem 3.7 of \cite{MS} we obtain the following result.

\begin{theorem}\label{th4}
The generating function  for the number of  weak ascent sequences of length $n$ that avoid  $102$ satisfies $f=1+x(1-x)f^3$.
\end{theorem}

The  series expansion of the generating function $f$ is
\begin{align*}
1+ x + 2 x^2 + 6 x^3 + \bm{22} x^4 + 89 x^5 + 381 x^6 + 1694 x^7 + 7744 x^8 +O(x^9).
\end{align*}
The weak ascent sequences corresponding to the  bold coefficient in the above series are
\begin{multline*}
\A_4(102)=\{ 0000, 0001, 0002, 0003, 0010, 0011, 0012, 0013, 0020, 0021,\\ 0022, 0023, 0100, 0101, 0
  110, 0111, 0112, 0113, 0120, 0121, 0122, 0123\}.
\end{multline*}

The sequence $\wa_{102}(n)$ coincides  with  the sequence  \seqnum{A200753} and
$$\wa_{102}(n)=\sum_{k=0}^{\lfloor \frac n2 \rfloor} \frac{(-1)^k}{2(n - k) + 1}\binom{n-k}{k}\binom{3(n - k)}{n - k},  \text{ for } n\geq 1.$$

\subsection{Pattern 012}
The pattern $012$ also leads to a well-known combinatorial number sequence. The weak ascent sequences avoiding pattern $012$ are enumerated by the  odd  indexed Fibonacci numbers.  Weak ascent sequences of length $4$ avoiding the pattern $012$ are listed below
\begin{align*}
\A_4(012)=\{ 0000, 0001, 0002, 0003, 0010, 0011, 0020, 0021, 0022, 0100, 0101, 0110, 0111\}.
\end{align*}
We find the generating function of the sequence  $\wa_{012}(n)$ by analyzing the maximal prefix of zeros of the elements in $\A_n(012)$.

Let $F_{\tau;a_1a_2\cdots a_s}(x)$ be the generating function for the number of weak ascent sequences $\pi=\pi_1\pi_2\cdots\pi_n$ that avoid $\tau$  and have prefix $a_1a_2\cdots a_s$. Note that $F_{\tau}(x)=F_{\tau;\epsilon}(x)$  is the generating function for the number of weak ascent sequences that avoid $\tau$, where $\epsilon$ is the empty word. By the definitions we have that  $F_{\tau}(x)=1+F_{\tau;0}(x)$.

\begin{theorem}\label{th1}
The generating function for the number of weak ascent sequences of length $n$ that avoid  $012$ is given by
$$F_{012}(x) = \frac{1-2x}{1-3x+x^2}.$$
\end{theorem}
\begin{proof}
Let $x$ be a weak ascent sequence in $\A(012)$ and have prefix $0^m$ for some  integer $m\geq 0$.  Then, we have the decomposition $x=0^mx'$, where $x'\in \A(012)$. If $x'=\epsilon$, then the generating function for this case is $x^m$. If $x'\neq \epsilon$, then the letter $x_{m+1}$ of $x$ satisfies $0\leq x_{m+1}\leq \wasc(0^m)+1=(m-1)+1=m$. This implies the functional equation:
\begin{align*}
F_{012;0^m}(x)&=x^m+F_{012;0^{m+1}}(x)+\sum_{j=1}^mF_{012;0^mj}(x).
\end{align*}
From a similar argument we have the equality
\begin{align*}
F_{012;0^mj}(x)&=x^{m+1}+F_{012;0^mj0}(x)+\sum_{i=1}^jF_{012;0^mji}(x)\\
&=x^{m+1}+xF_{012;0^mj}(x)+x\sum_{i=1}^jF_{012;0^mi}(x).
\end{align*}
By induction on $j=1,2,\ldots,m$, we obtain the expression
\begin{align*}
F_{012;0^mj}(x)&=\frac{x^{m+1}(1-x)^{j-1}}{(1-2x)^j}.
\end{align*}
Thus,
\begin{align*}
F_{012;0^m}(x)&=x^m+F_{012;0^{m+1}}(x)+\sum_{j=1}^m\frac{x^{m+1}(1-x)^{j-1}}{(1-2x)^j}=\frac{x^m(1-x)^m}{(1-2x)^m}+F_{012;0^{m+1}}(x),
\end{align*}
which leads to
\[F_{012}(x)=\sum_{m\geq0}\frac{x^m(1-x)^m}{(1-2x)^m}=\frac{1-2x}{1-3x+x^2}. \qedhere\]
\end{proof}

The  series expansion of the generating function $F_{012}(x)$ is
\begin{align*}
1+x + 2 x^2 + 5 x^3 + \bm{13} x^4 + 34 x^5 + 89 x^6 + 233 x^7 + 610 x^8 +O(x^9).
\end{align*}

The weak ascent sequences corresponding to the bold coefficient were given above.

We note that the generating function $F_{012}(x)$ coincides with the generating function of the sequence $F_{2n-1}$, where $F_n$ denotes the $n$-th Fibonacci number and we also give a combinatorial proof for this fact.

\subsubsection{Combinatorial proof}

First, we mention that $012$-avoiding weak ascent sequences are $012$-avoiding inversion sequences since a $012$-avoiding inversion sequence always has  the weak ascent property. So $\A_n(012)=\I_n(012)$ for all positive integer $n$. \cite{CMSW} considered already these sequences and proved that the numbers $\wa_{012}(n)$ satisfy the recurrence
$\wa_{012}(n) = 3\wa_{012}(n-1)-\wa_{012}(n-2)$, for $n\geq 3$, with initial conditions $\wa_{012}(1)=1$ and $\wa_{012}(2)=2$. They point out that the so-called Boolean permutations are also equinumerous and because of the connections to Coxeter groups it would be nice to find a simple bijection. Here, we present a bijection between $012$-avoiding weak ascent sequences (or inversion sequences) and a set that is well known to be enumerated by the odd-indexed Fibonacci numbers. The motivation of our bijection is the characterization given in Observation 1 in  \cite{CMSW}. 
\begin{observation}[\cite{CMSW}]
	Inversion sequences that avoid the pattern $012$ are those whose positive elements form a weakly decreasing sequence. 
\end{observation}
On the other hand, several known classical combinatorial objects are counted by the Fibonacci numbers $F_{2n-1}$. For instance, the number of matchings of a path on $2n$ vertices, the number of coverings of a $1\times 2n$ rectangle by $1\times 1$ squares and $1\times 2$ rectangles, the number of compositions of $2n$ by the numbers of $1$ and $2$. We formulate our result in terms of $1-2$ compositions.

First, we define three operations and show that each $012$-avoiding weak ascent sequence is uniquely determined by a sequence of these operations. Hence, there is a one-to-one correspondence between sequences of such operations and $012$ avoiding weak ascent sequences. Then, we show that an operation-sequence can be simply decoded by a $1-2$ composition.

{The operations $O$, $I$, and $E_{\ell}$, where $\ell$ is an integer $\geq 2$, are defined as follows:
\begin{itemize}
\item[$O$:] write a $0$ at the end of the sequence.
\item[$I$:] write a $1$ at the end of the sequence.
\item[$E_{\ell}$:] write at the beginning of the sequence $\ell-1$ $0$'s, increase each positive element by $\ell-1$, and write $\ell$ at the end of the sequence ($\ell\geq 2$).
\end{itemize}
It is clear that if $e$ is a $012$-avoiding weak ascent sequence, then $O(e)$, $I(e)$, and $E_{\ell}(e)$ are also $012$-avoiding weak ascent sequences. Moreover, if $e$ has length $n$, then $O(e)$ and $I(e)$ have length $n+1$, and $E_{\ell}(e)$ has length $n+\ell$. On the other hand, the operation-sequence which decodes a given $012$ avoiding inversion sequence is uniquely defined.

\begin{example}
We demonstrate how the operation-sequence $E_3OOIE_2OOIIE_4$ generates the corresponding weak ascent sequence. Let $\epsilon$ denote the empty weak ascent sequence. Then
\begin{align*}
E_4(\epsilon)\quad&\rightarrow\quad 0004\\
IE_4(\epsilon)\quad&\rightarrow\quad 00041\\
IIE_4(\epsilon)\quad&\rightarrow\quad 000411\\
OIIE_4(\epsilon)\quad&\rightarrow\quad 0004110\\
OOIIE_4(\epsilon)\quad&\rightarrow\quad 00041100\\
E_2OOIIE_4(\epsilon)\quad&\rightarrow\quad 0000522002\\
OE_2OOIIE_4(\epsilon)\quad&\rightarrow\quad 00005220020\\
IOE_2OOIIE_4(\epsilon)\quad&\rightarrow\quad 000052200201\\
E_3IOE_2OOIIE_4(\epsilon)\quad&\rightarrow\quad 000000744004033.
\end{align*}
Finally, we have to add a $0$ at the beginning of the so-created sequence.
\end{example}

Now we encode the operation-sequence by a composition of $1$'s and $2$'s. Write for a  $O$ two $1$'s: $\texttt{1}\texttt{1}$, for a $I$, write a $\texttt{2}$, and for an operation $E_{\ell}$ write $\ell-1$ $2$'s between two $1$'s, i.e.,
$\texttt{12}^{\ell-1}\texttt{1}$. Note that this encoding is unique, hence the map is invertible. 
Our example, $E_3IOE_2OOIIE_4$, translates to the composition
\begin{align*}\texttt{1221}\cdot \texttt{2}\cdot\texttt{11}\cdot\texttt{121} \cdot\texttt{11}\cdot\texttt{11}\cdot\texttt{2}\cdot \texttt{2}\cdot\texttt{12221}.\end{align*}

The above bijection is based on the recursive structure of the elements in $\A_n(012)$, which we could formulate as follows. (This bijection was pointed out by one of the referees.) Let $x\in \A_n(012)$. If the last letter of $x$ is $0$ or $1$, remove this last letter, obtaining a sequence of the set $\A_{n-1}(012)$. If the last letter of $x$ is $k>1$, remove that last letter, remove the first $k-1$ zeros from $x$, and subtract $(k-1)$ from each remaining positive letter, to obtain such a sequence of the set $\A_{n-k}(012)$.

\subsection{Pattern 021}
The generating function for the case $021$ can be derived using the previous technique, however, the derivation and the result are more complicated.
\begin{theorem}\label{th2}
The generating function $F_{021}(x)$ for the number of weak ascent sequences of length $n$ that avoid  $021$ satisfies
$F_{021}(x)=1+xF_{021}(x)/(1-xF_{021}(x)^2)$, which implies
$$F_{021}(x)-1=-\frac{xv_0^2-v_0+x+1}{x},$$
where
$$v_0=\frac{2\sqrt{1-4x-2x^2}}{3x}\cos\left(-\frac{1}{3}\arccos\left(\frac{7x^3-6x^2+12x-2}{2(1-4x-2x^2)^{3/2}}\right)+\frac{\pi}{3}\right)+\frac{x+1}{3x}.$$
\end{theorem}
\begin{proof}
By the definitions, we have
\begin{align}
F_{021;0^m}(x)&=x^m+F_{021;0^{m+1}}(x)+\sum_{j=1}^mF_{021;0^mj}(x)\nonumber\\
&=x^m+F_{021;0^{m+1}}(x)+\sum_{j=1}^mx^{j-1}F_{021;0^{m+1-j}1}(x),\label{eq021a}
\end{align}
where we used here that $F_{021;0^mj}(x)=x^{j-1}F_{021;0^{m+1-j}1}(x)$. To see this fact consider a weak ascent sequence $\pi=0^mj\pi'$ of length $n$. Since $\pi$ avoids $021$,   $\pi'$ does not contain any of the entries $1,2,\ldots,j-1$. In other words, each positive entry is at least $j$. Thus, $\pi$ avoids $021$ if and only if $0^{m+1-j}1\pi''$ avoids $021$, where $\pi''$ is obtained from $\pi'$ by decreasing each positive entry  of $\pi'$ by $j-1$. Clearly, $0^{m+1-j}1\pi''$ is a weak ascent sequence.

The generating function $F_{021;0^m1}(x)$ satisfies
\begin{align*}
F_{021;0^m1}(x)&=x^{m+1}+F_{021;0^m10}(x)+\sum_{j=1}^{m+1}F_{021;0^m1j}(x).
\end{align*}
Clearly, $F_{021;0^m11}(x)=F_{021;0^{m+1}1}(x)$. As before, we have that
$F_{021;0^m1j}(x)=x^{j-1}F_{021;0^{m+2-j}1}(x)$, for all $j=2,3,\ldots,m$. So,
\begin{align}
F_{021;0^m1}(x)&=x^{m+1}+F_{021;0^m10}(x)+\sum_{j=1}^{m+1}x^{j-1}F_{021;0^{m+2-j}1}(x).\label{eq021b}
\end{align}
Similarly, we have
\begin{align}
F_{021;0^m10}(x)&=x^{m+2}+F_{021;0^m100}(x)+\sum_{j=1}^{m+1}F_{021;0^m10j}(x)\nonumber\\
&=x^{m+2}+F_{021;0^{m+1}10}(x)+\sum_{j=1}^{m+1}x^jF_{021;0^{m+2-j}1}(x).\label{eq021c}
\end{align}
In order to solve these recurrences, we define the generating functions $A(v)=\sum_{m\geq1}F_{021;0^m}(x)v^{m-1}$,   $B(v)=\sum_{m\geq1}F_{021;0^m1}(x)v^{m-1}$, and $C(v)=\sum_{m\geq1}F_{021;0^m10}(x)v^{m-1}$. Clearly, $A(0)=F_{021;0}(x)=F_{021}(x)-1$. So, our next step is to find an explicit formula for the generating function $A(0)$. Thus, by multiplying \eqref{eq021a}-\eqref{eq021c} by $v^{m-1}$ and summing over $m\geq1$, we obtain
\begin{align}
A(v)&=\frac{x}{1-xv}+\frac{1}{v}(A(v)-A(0))+\frac{1}{1-xv}B(v),\label{eq021aa1}\\
B(v)&=\frac{x^2}{1-xv}+C(v)+\frac{1}{v(1-xv)}B(v)-\frac{1}{v}B(0),\label{eq021aa2}\\
C(v)&=\frac{x^3}{1-xv}+\frac{1}{v}(C(v)-C(0))+\frac{x}{v(1-xv)}B(v)-\frac{x}{v}B(0).\label{eq021aa3}
\end{align}
By \eqref{eq021aa1}-\eqref{eq021aa2}, we have
\begin{multline}
(v^3x-v^2x-v^2+vx+2v-1)B(v)\\
=(v^2x^2+v^2x-2vx-v+1)B(0)+(v^2x-v)C(0)+v^2x^3+v^2x^2-vx^2.\label{eq021aa4}
\end{multline}
Let $v_0,v_1,v_2$ be the roots of the polynomial equation $v^3x-v^2x-v^2+vx+2v-1=0$, namely\footnote{See the following link to find the roots in terms of trigonometric functions: \url{https://en.wikibooks.org/wiki/Trigonometry/The_solution_of_cubic_equations}},
$$v_j=\frac{2\sqrt{1-4x-2x^2}}{3x}\cos\left(-\frac{1}{3}\arccos\left(\frac{7x^3-6x^2+12x-2}{2(1-4x-2x^2)^{3/2}}\right)+\frac{\pi(2j+1)}{3}\right)+\frac{x+1}{3x}.$$
By substituting $v=v_1$ and $v=v_2$ into \eqref{eq021aa4}, we obtain the following system of equations:
\begin{align*}
&(v_1^2x^2+v_1^2x-2v_1x-v_1+1)B(0)+(v_1^2x-v_1)C(0)+v_1^2x^3+v_1^2x^2-v_1x^2=0,\\
&(v_2^2x^2+v_2^2x-2v_2x-v_2+1)B(0)+(v_2^2x-v_2)C(0)+v_2^2x^3+v_2^2x^2-v_2x^2=0.
\end{align*}
By solving this system for $B(0)$ and $C(0)$, we obtain that
$$B(0)=\frac{v_1v_2x^2}{(v_2x-1)(v_1x-1)} \text{\quad and \quad } C(0)= -\frac{(v_2x+v_2-1)(v_1x+v_1-1)x^2}{(v_2x-1)(v_1x-1)},$$
where we simplified the expressions by using the fact that $v=v_1$ and $v=v_2$ satisfy $v^3x-v^2x-v^2+vx+2v-1=0$. 
Since $v_0,v_1,v_2$ are solutions of $v^3x-v^2x-v^2+vx+2v-1=0$, so  $v_0+v_1+v_2=(1+x)/x$ and $v_0v_1v_2=1/x$. Thus,
$$B(0)=\frac{1}{v_0^2-v_0+1} \text{\quad and \quad } C(0)=-\frac{v_0^2x(x+1)-v_0(x^2+x+1)+x^2+2x+1}{v_0^2-v_0+1}.$$
Now, by substituting these expressions into \eqref{eq021aa2}-\eqref{eq021aa3} and solving for $B(v)$ and $C(v)$, we obtain
$$B(v)=-\frac{vxv_0^2+v(xv-x-1)v_0+1}{(v_0^2-v_0+1)(v^3x-v^2(x+1)+v(x+2)-1)}.$$
Now, by \eqref{eq021aa1} with $v=1$, we have
$$A(0)=\frac{x}{1-x}+\frac{1}{1-x}B(1),$$
which, after simplification, implies that
$$A(0)=-\frac{xv_0^2-v_0+x+1}{x}=F_{021}(x)-1.$$
Note that $v_0$ is a root of the polynomial $v^3x-v^2x-v^2+vx+2v-1$, so $F_{021}(x)$ satisfies
\[F_{021}(x)=1+\frac{xF_{021}(x)}{1-xF_{021}(x)^2}. \qedhere\]
\end{proof}

The  series expansion of the generating function $F_{021}(x)$ is
\begin{align*}
1+x + 2 x^2 + 6 x^3 + \bm{21} x^4 + 80 x^5 + 322 x^6 + 1347 x^7 + 5798 x^8  +O(x^9).
\end{align*}
The weak ascent sequences corresponding to the  bold coefficient in the above series are
\begin{multline*}
\A_4(021)=\{ 0000, 0001, 0002, 0003, 0010, 0011, 0012, 0013, 0020, 0022, 0023, \\
0100, 0101, 0102, 0110, 0111, 0112, 0113, 0120, 0122, 0123 \}.
\end{multline*}

The sequence $\wa_{021}(n)$ coincides  with  the sequence \seqnum{A106228}, and
$$\wa_{021}(n)=\sum_{k=0}^{n-1} \frac{1}{n}\binom{2 n - 2 k - 2}{n - k - 1}\binom{n + k - 1}{n - 1}, \text{ for } n\geq 1.$$
It is known that node-labeled plane trees are enumerated by these numbers. We present a bijection between weak ascent sequences avoiding $021$ and node-labeled plane trees. However, in order to have a more insightful picture, we consider a special matrix class that is in bijection with weak ascent sequences.
\subsubsection{Matrices} The structure of a weak ascent sequence is nicely recognizable in a certain class of upper triangular $01$ matrices, $\newmatrices$ defined in \cite{BCD}. We recall the definition and some important properties.

\begin{definition}\label{matrixdefn}
  Let $\newmatrices_n$ be the set of upper triangular square
  $01$-matrices $A$ that satisfy the following properties:
  \begin{enumerate}
  \item[(a)] There are $n$ 1's in $A$.
  \item[(b)] There is at least one 1 in every column of $A$.
  \item[(c)] For every pair of adjacent columns, the top-most 1 in the left column is weakly above the bottom-most 1 in the right column.
\end{enumerate}
\end{definition}
There is a one-to-one correspondence between the set $\newmatrices_n$ and the set of weak ascent sequences of length $n$ as it is shown in \cite{BCD}. Given a weak ascent sequence, $w = 00211015$, the columns of the corresponding matrix $M\in \newmatrices_n$ are constructed according to the decreasing subsequences (\emph{decreasing runs}) of $w$. The decomposition into decreasing subsequences is unique, in our example $w = 0|0|21|10|1|5$. For each decreasing run there is a column with $1$ entries in the rows given in the particular subsequences, i.e., the $1$ entries of the matrix are in positions $(i,j)$ if the $j$th decreasing run contains the value $i-1$.

Forbidding  a certain pattern in weak ascent sequences affects the structure of the corresponding matrix. In certain cases, this structure is nice and gives an illuminating insight. This is the case for the pattern $021$ which we describe now.

Let $\newmatrices(021)$ denote the set of matrices corresponding to weak ascent sequences avoiding the pattern $021$. Recall that weak ascent sequences that avoid $021$ are those whose positive elements form a non-decreasing sequence. As a consequence, in a matrix $M\in \newmatrices(021)$ the $1$ entries are of two types. The $1$'s corresponding to the  $0$'s of the weak ascent sequence are all in the top row. The other $1$'s corresponding to the positive elements of the weak ascent sequence form a weakly increasing staircase path.

This characterization is sufficient. All $M\in \newmatrices$ whose $1$ entries are either in the top row or form a weakly increasing chain are exactly the matrices of $\newmatrices(021)$.

Note that our matrices are upper triangular, so each matrix can be associated with a path having all steps above the diagonal, in other words with a Dyck path, with eventually some extra $1$'s in the first row.
In particular, matrices with dimension $n$ that do not contain any $1$'s in the first row (besides the $1$ entry in the first column that is contained in all matrices) can be associated to Dyck paths. Hence, the subclass of $n\times n$ matrices in $\newmatrices_n(021)$ is enumerated by the $n$-th Catalan number. (Note that in the set $\newmatrices_n(021)$ there are matrices of different dimensions.)

\subsubsection{Bijection with node-labeled plane trees}
We describe now a bijection between $021$ avoiding weak ascent sequences and node-labeled plane trees on $n$ vertices. A \emph{plane tree} is an oriented tree in which the children of a vertex are ordered.

Plane trees are enumerated by Catalan numbers and their combinatorics is very rich and popular, see for example  \cite{SC}. We recall some facts and decodings that are useful for our purposes. See for more details for instance \cite{DZ}.

A plane tree can be decoded by a sequence $(a_1,\ldots, a_n)$ of non-negative integers satisfying $a_1+\cdots + a_n\geq i$ and $\sum_{j=0}^n a_j = n$. See Exercise 86, page 33 in \cite{SC}.	
 We obtain the corresponding sequence to a plane tree by traversing the tree in preorder (visiting the root, then traversing its subtrees from left to right) and recording the outdegree of the node we just visited (the first time). For $n=3$ the five sequences are
\begin{align*}
(3,0,0)\quad (2,1,0) \quad (2,0,1) \quad (1,2,0) \quad (1,1,1).
\end{align*}
We can associate a Dyck path to the sequence $(a_1,\ldots, a_n)$ by the following rule:
\begin{align*}
	H^{a_1}DH^{a_2}D\cdots H^{a_{n}}D,
\end{align*}
where $H$ denotes a horizontal step and $D$ a down step.
We call a vertex in the plane tree that is not the root a \emph{node} if its outdegree is greater than zero.  We label the node $v$ of the tree with an integer $\ell_v$ such that $0\leq \ell_v\leq \mbox{odeg}(v)-1$, where $\mbox{odeg}(v)$ denotes the outdegree of $v$. A plane tree with such a labeling of all nodes is called \emph{node-labeled plane tree}. Note that we slightly modified  the definition and terminology of David Callan, see sequence  \seqnum{A106228} in  \cite{OEIS}.

To each node-labeled plane tree we can associate a sequence $(a_1,a_2[\ell_2]\ldots, a_n[\ell_n])$ such that \linebreak $\sum_{j = 0}^i a_i\geq i$, $\sum_{j=0}^n a_j = n$, and $0\leq \ell_i< a_i$ for all $1\leq i \leq n$. We call such a sequence an \emph{augmented sequence}. (Note that $\ell_1$ is not defined since there is no label associated to the root.) The $6$ augmented sequences of length $3$ are listed below.
\begin{align*}
		(3,0,0)\quad (2,1[0],0)\quad (2,0,1[0]) \quad (1,2[0],0) \quad (1,2[1],0) \quad (1,1[0],1[0]).
\end{align*}
In the plane tree we can think of the label as a mark on one of the descendants of the node, whereas in the Dyck path as a mark on one of the horizontal steps. More precisely, we define a \emph{marked Dyck path} associated to an augmented sequence $(a_1,a_2[\ell_2],\ldots, a_n[\ell_n])$ as follows:
\begin{align*}
	H^{a_1}DH^{a_2-\ell_2}h^{\ell_2}D\cdots H^{a_{n}- \ell_n}h^{\ell_n}D.
\end{align*}
For instance, if the augmented sequence is $(1,3[1],2[1],0,0,4[3],0,0, 1[0], 0, 1[0], 0)$, the corresponding marked Dyck path is $HDHHhDHhDDDHhhhDDDHDDHDD$.

We associate now to a marked Dyck path a $01$ matrix. If all $\ell_i = 0$, then write for each horizontal step a $1$ in the row of its height. More precisely, in row $r$, where $r$ is given by the number of down steps standing to the left of the particular horizontal step. This way we obtain the $n\times n$ matrices of the set $\A_n(021)$. We illustrate this in the Example \ref{ex_matrix_Dyck_path}.
\begin{example}\label{ex_matrix_Dyck_path}
Let's consider the case $n=4$.  Of the 21 matrices in $\newmatrices_n(021)$, fourteen have a single 1 in the first row. These $14$ cases correspond to those labeled plane trees where all $\ell_i=0$. The corresponding marked Dyck paths do not contain any $h$. For instance, the sequence $(2,1[0],1[0])$ corresponds to  the Dyck path $HHDHDHDD$:
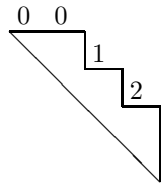
\begin{center}
\setlength{\unitlength}{0.5cm}
\begin{picture}(6,6)
\put(1,5){\line(1,-1){4}}
\put(1,5){\line(1,0){2}}
\put(3,4){\line(1,0){1}}
\put(4,3){\line(1,0){1}}
\put(3,5){\line(0,-1){1}}
\put(4,4){\line(0,-1){1}}
\put(5,3){\line(0,-1){2}}
\put(1.2,5.2){$0$}\put(2.2,5.2){$0$}\put(3.2,4.2){$1$}
\put(4.2,3.2){$2$}
\end{picture}
\end{center}
The matrix is given by putting a $1$ entry for each horizontal step at its height
$
\left[\begin{smallmatrix}
1&1&0&0\\
0&0&1&0\\
0&0&0&1\\
0&0&0&0
\end{smallmatrix}\right]
$.
The weak ascent sequence is given by the heights of the horizontal steps $w = 0012$.
\end{example}
We now define the map also for  cases when there is at least one positive $\ell_i$. To a horizontal step, $h$, a $1$ entry will correspond in the first row of the matrix. If $\ell_i=1$, we apply the following rule when visiting the horizontal step $h$. If it would be a $H$ we would write a $1$ entry in the next column of the same row, but as it is a $h$ we write a $1$ entry in the first row of the same column. At the same time we ``delete" the $D$ followed by the $h$ from the path implying actually that the heights of the next horizontal steps are reduced by one.
For instance, $(1,2[1],1[0],0)$ transforms to the marked Dyck path  $HDHhDHDD$, the corresponding matrix is
	$\left[\begin{smallmatrix}1&1&0\\0&1&1\\0&0&0\end{smallmatrix}\right]$ and the weak ascent sequence is $0101$.

In the case of arbitrary $\ell_i$ we write $\ell_i$ $1$'s in the first row in consecutive columns. We delete the $D$ (coming after the $h$'s) as well, reducing with that the height of all the following $1$'s by one.
\begin{example}
For $n=4$ there are $7$ node-labeled plane trees with at least one positive label.
Here we give the list of them, the corresponding augmented sequences, and the corresponding matrices (in the same order).
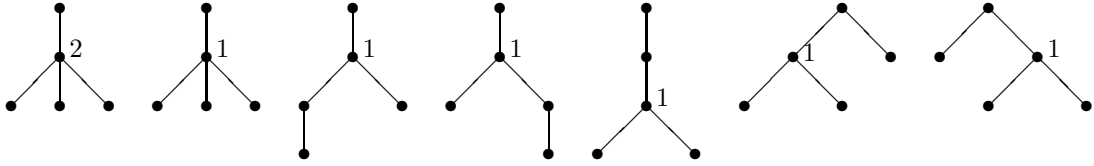
\begin{center}
\setlength{\unitlength}{0.65cm}
\begin{picture}(22,3)
\put(0,1){\circle*{0.2}}\put(1,1){\circle*{0.2}}\put(2,1){\circle*{0.2}}\put(1,2){\circle*{0.2}}\put(1,3){\circle*{0.2}}
\put(1.2,2){$2$}
\put(1,1){\line(0,1){2}}\put(0,1){\line(1,1){1}}\put(1,2){\line(1,-1){1}}
\put(3,1){\circle*{0.2}}\put(4,1){\circle*{0.2}}\put(5,1){\circle*{0.2}}\put(4,2){\circle*{0.2}}\put(4,3){\circle*{0.2}}
\put(4.2,2){$1$}
\put(4,1){\line(0,1){2}}\put(3,1){\line(1,1){1}}\put(4,2){\line(1,-1){1}}
\put(6,0){\circle*{0.2}}\put(6,1){\circle*{0.2}}\put(8,1){\circle*{0.2}}\put(7,2){\circle*{0.2}}\put(7,3){\circle*{0.2}}
\put(7.2,2){$1$}
\put(6,0){\line(0,1){1}}\put(6,1){\line(1,1){1}}\put(7,2){\line(0,1){1}}\put(7,2){\line(1,-1){1}}
\put(9,1){\circle*{0.2}}\put(10,2){\circle*{0.2}}\put(10,3){\circle*{0.2}}\put(11,0){\circle*{0.2}}\put(11,1){\circle*{0.2}}
\put(10.2,2){$1$}
\put(9,1){\line(1,1){1}}\put(10,2){\line(0,1){1}}\put(10,2){\line(1,-1){1}}\put(11,1){\line(0,-1){1}}
\put(12,0){\circle*{0.2}}\put(14,0){\circle*{0.2}}\put(13,1){\circle*{0.2}}\put(13,2){\circle*{0.2}}\put(13,3){\circle*{0.2}}
\put(13.2,1){$1$}
\put(12,0){\line(1,1){1}}\put(13,1){\line(0,1){2}}\put(13,1){\line(1,-1){1}}
\put(15,1){\circle*{0.2}}\put(16,2){\circle*{0.2}}\put(17,1){\circle*{0.2}}\put(17,3){\circle*{0.2}}\put(18,2){\circle*{0.2}}
\put(16.2,1.9){$1$}
\put(15,1){\line(1,1){2}}\put(16,2){\line(1,-1){1}}\put(17,3){\line(1,-1){1}}
\put(19,2){\circle*{0.2}}\put(20,3){\circle*{0.2}}\put(20,1){\circle*{0.2}}\put(21,2){\circle*{0.2}}\put(22,1){\circle*{0.2}}
\put(21.2,2){$1$}
\put(19,2){\line(1,1){1}}\put(20,3){\line(1,-1){2}}\put(20,1){\line(1,1){1}}
\end{picture}
\end{center}
\begin{align*}
&(1, 3[2],0,0) \quad (1,3[1],0,0)  \quad (1,2[1],1[0],0) \quad (1,2[1],0,1[0])\\
&(1, 1[0],2[1],0) \quad  (2, 2[1],0,0) \quad (2,0,2[1],0).
\end{align*}
\begin{align*}
&\left[\begin{matrix}1&1&1\\0&1&0\\0&0&0\end{matrix}\right]\quad \left[\begin{matrix}1&0&1\\0&1&1\\0&0&0\end{matrix}\right]\quad
\left[\begin{matrix}1&1&0\\0&1&1\\0&0&0\end{matrix}\right]\quad
\left[\begin{matrix}1&1&0\\0&1&0\\0&0&1\end{matrix}\right]\quad
\left[\begin{matrix}1&0&1\\0&1&0\\0&0&1\end{matrix}\right]\quad
\left[\begin{matrix}1&1&1\\0&0&1\\0&0&0\end{matrix}\right]\quad
\left[\begin{matrix}1&1&1\\0&0&0\\0&0&1\end{matrix}\right].
\end{align*}
\end{example}

\subsection{Pattern 011}
The main observation for the structure of $011$-avoiding weak ascent sequences is that each positive number can occur only once. We introduce some notation and terminology. Let $w=w_1\cdots w_n$ be a weak ascent sequence that avoids the pattern $011$.  If $w_{j-1}\leq w_j$ for a $j$, we denote by $q(w_j)$ the value that is newly available for $w_{j+1}$, $q(w_j) = \wasc(w_1w_2\cdots w_{j})+1$. We say that $w_j$ is a \emph{candidate} and $q(w_j)$ is its \emph{follower}. Further, let $s(w_1\cdots w_j)$ denote the set of all available values for $w_{j+1}$.

From the set $s_{j-1}$ we get the set $s_j$ depending on the value of $w_j$. Three cases are possible:
\begin{itemize}
	\item[A.] If $w_{j-1}=w_j=0$, then $s_j=s_{j-1}\cup\{q(w_j)\}$.
	\item[B.] If $w_{j-1}<w_j$, then  $s_j=s_{j-1}\cup\{q(w_j)\}\setminus\{w_{j}\}$.
	\item[C.] If $w_{j-1}>w_j$, then if $w_j=0$, we have $s_j=s_{j-1}$, and if $w_j>0$, we have $s_j=s_{j-1}\setminus\{w_j\}$.
\end{itemize}

If the follower $q(w_j)$ is also a candidate, meaning that it occurs in the weak ascent sequence and it is greater than its left neighbour (note that it can not be equal to its left neighbour) then we can consider also its follower, $q(q(w_j))$. We can continue this procedure, obtaining a sequence of values $w_j$, $q(w_j)$, $\ldots$, $q(\cdots q(w_j))$ until the last is not a candidate. This means that it is smaller than its left neighbour in the weak ascent sequence or it does not occur in the weak ascent sequence as an element. Hence, such a sequence is maximal if the first, starting candidate is zero and the last element is not a candidate. Moreover, each $0$ element in the weak ascent sequence determines such a maximal sequence and each positive number from the set $\{1,2,\ldots, q(w_n)\}$ arises in exactly one such maximal sequence. Hence, a set partition of $\{1,2,\ldots, q(w_n)\}$ uniquely defines a $011$-avoiding weak ascent sequence $w\in \A_n(011)$. More precisely, let $i_1$, $i_2$, $\ldots$, $i_m$ denote the positions of the zeros in the weak ascent sequence. For each $1\leq k\leq m$, $w_{i_k}=0$, define the set $S_k$ containing the elements of the maximal sequence of the followers:
\begin{align*}
	S_k = \{q(w_{i_k}), \ldots, q(q(\cdots(q(w_{i_k}))))\}.
\end{align*}

We illustrate the above definition by an example.
\begin{example}
	Let us consider the weak ascent sequence $w=0001300524$. We added in the following table the value $q(w_j)$ above the element $w_j$.
\begin{center}
\begin{tabular}{ccccccccccc}
$q(w_j)$&\textbf{1}&2&3&\textbf{4}&5&&6&7&&\textbf{8}\\
$w$&\textbf{0}&0&0&\textbf{1}&3&0&0&5&2&\textbf{4}
\end{tabular}
\end{center}
There are $5$ zeros in the weak ascent sequence. Starting with the first zero, we have $q(w_{i_1})=q(w_1)=q(0)=1$, $q(1)=4$, $q(4)=8$. Hence, the maximal sequence is $0-1-4-8$, and $S_1= \{1,4,8\}$. The other sets are $S_2=\{2\}$, $S_3=\{3,5,7\}$, $S_4 = \emptyset$, $S_5 = \{6\}$.
\end{example}
After understanding the structure of this class, we derive its size. First, we provide a finer counting according to the number of zeros. 

Let $\mathcal{B}(n,m)$ denote the set of weak ascent sequences avoiding the pattern $011$ of length $n$ with $m$ zeros and let $B(n,m)=|\mathcal{B}(n,m)|$. The initial condition is $B(n,1) = 1$.
We have the following proposition.
\begin{proposition} For all $1\leq m\leq n$,
\begin{align*}
B(n,m) = \sum_{k=m}^n B(k-1,m-1)\binom{n-k+m-1}{n-k}.
\end{align*}
\end{proposition}
\begin{proof}
Let $w\in \mathcal{B}(n,m)$ and $k$ denote the position of the last zero. Then $w_1\cdots w_{k-1}$ is in the set $\mathcal{B}(k-1,m-1)$. The values $w_i$, $i>k$ are all positive numbers. Each positive number that occurs is contained in one of the sets $S_1,\ldots, S_{m-1}$ or in the set that starts with the last zero at the $m$th position, $S_m$. The number of ways to place $n-k$ elements into $m$ sets (possibly with repetition) is counted by the binomial coefficient $\binom{n-k+m-1}{n-k}$.
\end{proof}

 \cite{Ruskey} showed that the above recurrence is satisfied by the triangular binomial coefficient, $\binom{n-1+\binom{m}{2}}{n-m}$, which equals $1$ when $m = 1$.
\begin{theorem}
The number of weak ascent sequences of length $n$ that avoid the pattern $011$ is given by
$$\wa_{011}(n)= \sum_{m=1}^{n}\binom{n-1+\binom{m}{2}}{n-m}.$$
\end{theorem}
In  \seqnum{A098569} we find several combinatorial interpretations of these numbers. Interestingly, though some of the interpretations are integer sequences, none of these coincide with weak ascent sequences that avoid the pattern $011$. The corresponding matrices $\newmatrices_n(011)$ are those in $\newmatrices_n$ that have at most one $1$ entry in each row except the first row.

\section{Open questions}

We leave as an open question the enumeration of the remaining cases of length 3. In Table \ref{table2} we give the first few values of the sequence $\wa_p(n)$ for $p\in \{000, 010, 100, 101, 110, 120, 201, 210\}$.

\begin{table}[ht]
\centering
\begin{tabular}{|c|c|} \hline
Pattern $p$ & Number of sequences avoiding $p$ \\ \hline \hline
000 &  1, 2, 5, 15, 51, 194, 809, 3667, 17892, 93408  \\ \hline
010 & 1, 2, 5, 15, 52, 202, 861, 3969, 19582, 102600   \\ \hline
100 &  1, 2, 6, 22, 94, 452, 2398, 13832, 85838, 568412   \\ \hline
101 &  1, 2, 6, 22, 93, 438, 2251, 12447, 73308, 456401   \\ \hline
110 &  1, 2, 6, 22, 92, 424, 2112, 11229, 63174, 373692   \\ \hline
120 &  1, 2, 6, 22, 91, 409, 1958, 9860, 51775, 281654   \\ \hline
201 & 1, 2, 6, 23, 104, 528, 2919, 17225, 107022, 693327   \\ \hline
210 & 1, 2, 6, 23, 104, 530, 2958, 17734, 112657, 750726   \\ \hline
\end{tabular}
\caption{Sequence $\wa_p(n)$ for $p\in \{000, 010, 100, 101, 110, 120, 201, 210\}$}\label{table2}
\end{table}

\begin{conjecture}
The sequence $\wa_{210}(n)$ coincides with the sequence  \seqnum{A117106}.
\end{conjecture}
Note that \seqnum{A117106}  enumerates permutations avoiding the vincular pattern $2$-$41$-$3$.

\acknowledgements
We thank the referees for their careful reading of the manuscript and valuable comments that led to improvements in the presentation of the paper greatly.

\bibliographystyle{abbrvnat}
\bibliography{biblo.bib}

\end{document}